\documentclass[11pt,reqno]{amsart}

\usepackage[a4paper,margin=30mm]{geometry}
\usepackage{amsmath,amssymb,mathtools,mathrsfs}
\usepackage{enumitem}
\usepackage{booktabs}
\usepackage{aliascnt}
\usepackage{microtype}
\usepackage{xcolor}
\usepackage[colorlinks=true,linkcolor=blue!55!black,citecolor=blue!55!black,urlcolor=blue!55!black]{hyperref}
\hypersetup{
  pdftitle={A Moment Proof of the Hollmann--Xiang Conjectures on Kloosterman Polynomials, Thirteenth Version},
  pdfauthor={Yi Li},
  pdfsubject={Classification and structural reduction of Kloosterman polynomials over binary finite fields},
  pdfkeywords={permutation polynomial, Kloosterman polynomial, linearized polynomial, trace fiber, moment method, spectral zero}
}

\allowdisplaybreaks
\setlist[enumerate]{itemsep=2pt,topsep=4pt}

\newtheorem{theorem}{Theorem}[section]
\newaliascnt{proposition}{theorem}
\newtheorem{proposition}[proposition]{Proposition}
\aliascntresetthe{proposition}
\newaliascnt{lemma}{theorem}
\newtheorem{lemma}[lemma]{Lemma}
\aliascntresetthe{lemma}
\newaliascnt{corollary}{theorem}
\newtheorem{corollary}[corollary]{Corollary}
\aliascntresetthe{corollary}
\newaliascnt{claim}{theorem}

\aliascntresetthe{claim}
\newaliascnt{conjecture}{theorem}
\newtheorem{conjecture}[conjecture]{Conjecture}
\aliascntresetthe{conjecture}
\theoremstyle{definition}
\newaliascnt{definition}{theorem}

\aliascntresetthe{definition}
\theoremstyle{remark}
\newaliascnt{remark}{theorem}

\aliascntresetthe{remark}
\usepackage[nameinlink,noabbrev]{cleveref}
\crefname{theorem}{theorem}{theorems}
\Crefname{theorem}{Theorem}{Theorems}
\crefname{proposition}{proposition}{propositions}
\Crefname{proposition}{Proposition}{Propositions}
\crefname{lemma}{lemma}{lemmas}
\Crefname{lemma}{Lemma}{Lemmas}
\crefname{corollary}{corollary}{corollaries}
\Crefname{corollary}{Corollary}{Corollaries}
\crefname{claim}{claim}{claims}
\Crefname{claim}{Claim}{Claims}
\crefname{conjecture}{conjecture}{conjectures}
\Crefname{conjecture}{Conjecture}{Conjectures}
\crefname{definition}{definition}{definitions}
\Crefname{definition}{Definition}{Definitions}
\crefname{remark}{remark}{remarks}
\Crefname{remark}{Remark}{Remarks}

\newcommand{\F}{\mathbb F}
\newcommand{\Tr}{\operatorname{Tr}}
\newcommand{\supp}{\operatorname{supp}}

\newcommand{\T}{\mathcal T_1}
\newcommand{\Zmod}{\mathbb Z/n\mathbb Z}
\newcommand{\R}{\mathcal R_n}
\newcommand{\E}{\mathcal E}
\newcommand{\ind}{\mathbf 1}

\title[The Hollmann--Xiang conjectures]{A Moment Proof of the Hollmann--Xiang Conjectures on Kloosterman Polynomials}
\author{Yi Li}
\address{School of Computer Science, Shanghai Jiao Tong University, Shanghai, China}
\email{yi.li@sjtu.edu.cn}
\date{}

\begin{document}

\begin{abstract}
Let $F=\F_{2^n}$ and let
\[
 \mathcal T_1=\{t\in F:\Tr_{F/\F_2}(t)=1\}.
\]
For an integer $u=\sum_{i=0}^{n-1}u_i2^i$ of even binary weight, put
\[
 L_u(X)=\sum_{i=0}^{n-1}u_iX^{2^i},
 \qquad
 L_{1,u}(t)=t+L_u(t^{-1}).
\]
Hollmann and Xiang conjectured that $L_{1,u}$ permutes $\mathcal T_1$ if and only if
$u\in\{0,3,6,10\}$.  We prove the necessity by a moment method; the sufficiency is their
Theorem~4.1.  The quadratic moment yields a cyclic autocorrelation identity, while the cubic
moment gives a two-variable group-algebra identity.  An elimination argument reduces the
coefficient word to a monomial or a complementary monomial, and explicit cubic and quartic
moments complete the classification.

We also prove Hollmann and Xiang's Conjecture~4.5.  Namely, if $n>3$ and
\[
 L_{c,d}(x)=L_c(x)+L_d(x^{-1})
\]
is a Kloosterman polynomial on $\F_{2^n}$, then its positive linearized part $L_c$ is
injective on both trace fibers.  This second result follows from relative quadratic and cubic
moment identities and a spectral-zero rigidity argument.  Consequently every Kloosterman
polynomial is obtained by composing one of $L_{1,0},L_{1,3},L_{1,6},L_{1,10}$ with a
linearized map that is injective on both trace fibers.
\end{abstract}

\maketitle

\section{Introduction}

Permutation polynomials and permutation rational functions over finite fields form a meeting
point of finite-field arithmetic, combinatorics, coding theory, and cryptography.  In
characteristic two, linearized polynomials are especially natural: they encode $\F_2$-linear
maps and interact directly with the absolute trace.  For the families studied here, this
interaction places the problem on the two affine trace fibers, where Frobenius exponents and
cyclic binary carries can be analyzed through power moments.

The problem originates in the work of Hollmann and Xiang on binary Kloosterman-sum
identities \cite{HollmannXiang2004}.  For integers $c,d$ with binary expansions
$c=\sum c_i2^i$ and $d=\sum d_i2^i$, they introduced
\[
 L_c(X)=\sum_i c_iX^{2^i},\qquad L_d(X)=\sum_i d_iX^{2^i},
\]
and
\[
 L_{c,d}(x)=L_c(x)+L_d(x^{-1}).
\]
When the binary weight of $d$ is even, they call $L_{c,d}$ a Kloosterman polynomial if it is
injective on the trace-one set.  They proved that $L_{1,3}$, $L_{1,6}$, and $L_{1,10}$ are
Kloosterman polynomials for every extension degree \cite[Theorem~4.1]{HollmannXiang2004},
and formulated the following classification \cite[Conjecture~4.4]{HollmannXiang2004}.

\begin{conjecture}[Hollmann--Xiang]\label{conj:HX}
Let $0\le d<2^n$ have even binary weight.  Then
\[
 L_{1,d}(t)=t+L_d(t^{-1})
\]
permutes the trace-one elements of $\F_{2^n}$ if and only if
\[
 d\in\{0,3,6,10\}.
\]
\end{conjecture}

Hollmann and Xiang also isolated a structural reduction for general Kloosterman polynomials.
For $e\in\F_2$, write
\[
 \mathcal T_e=\{x\in\F_{2^n}:\Tr(x)=e\},
\]
and let $\mathscr L_n$ be the collection of binary linearized polynomials that are injective
on both $\mathcal T_0$ and $\mathcal T_1$.  Their second conjecture is
\cite[Conjecture~4.5]{HollmannXiang2004}:

\begin{conjecture}[Hollmann--Xiang reduction]\label{conj:HX-general}
Let $n>3$.  If $L_{c,d}$ is a Kloosterman polynomial on $\F_{2^n}$, then
$L_c\in\mathscr L_n$.
\end{conjecture}

The Kloosterman-polynomial framework was developed further by Cao, Hollmann, and Xiang
\cite{CaoHollmannXiang2008}.  Equivalent formulations in terms of permutation polynomials
were investigated by Yuan, Ding, Wang, and Pieprzyk \cite{YuanDingWangPieprzyk2008} and by
Singh and Singh \cite{SinghSingh2019}.  More recently, Bartoli, Li, and Zhou
\cite{BartoliLiZhou2025} established the classification in an asymptotic range in which the
extension degree is sufficiently large relative to the support of the relevant linearized
polynomial.  The present paper proves both conjectures without a size hypothesis.

Let
\[
 F=\F_{2^n},\qquad \Tr=\Tr_{F/\F_2},\qquad
 \T=\{t\in F:\Tr(t)=1\}.
\]
For $u=\sum_{i=0}^{n-1}u_i2^i$, put
\[
 L_u(X)=\sum_{i=0}^{n-1}u_iX^{2^i}.
\]
Our first main result is the complete classification in \Cref{conj:HX}.

\begin{theorem}[Classification theorem]\label{thm:main}
Let $n\ge1$, and let $0\le u<2^n$ have even binary weight.  The following are equivalent:
\begin{enumerate}[label=\textup{(\roman*)}]
 \item the Kloosterman polynomial
 \[
  L_{1,u}(t)=t+L_u(t^{-1})
 \]
 permutes $\T$;
 \item $u\in\{0,3,6,10\}$, subject to $u<2^n$.
\end{enumerate}
\end{theorem}

The sufficiency in \Cref{thm:main} is exactly
\cite[Theorem~4.1]{HollmannXiang2004}, together with the trivial case $u=0$.  Sections
\ref{sec:powersums}--\ref{sec:completion} prove the necessity.  The proof is based entirely
on moments on the trace-one hyperplane.  The quadratic moment first places the coefficient
word in the unit group of a cyclic group algebra.  The cubic moment then yields a
 two-variable identity whose elimination forces the word to be a monomial or the complement
of a monomial.  Explicit cubic moments eliminate the complementary case; in the monomial
case another cubic moment and a quartic moment leave only $u=0,3,6,10$.

Our second main result proves \Cref{conj:HX-general}.

\begin{theorem}[General reduction theorem]\label{thm:general-reduction}
Let $n>3$, and let $c,d$ have binary expansions of length $n$, with $d$ of even binary
weight.  If
\[
 L_{c,d}(x)=L_c(x)+L_d(x^{-1})
\]
is a Kloosterman polynomial on $\F_{2^n}$, then $L_c\in\mathscr L_n$.
\end{theorem}

The proof of \Cref{thm:general-reduction} uses a relative version of the same moment method.
For the coefficient polynomials $C(Z)$ and $D(Z)$, set
\[
 B(Z)=D(Z)+ZC(Z).
\]
The quadratic moment gives
\[
 B(Z)B(Z^{-1})=C(Z)C(Z^{-1}),
\]
while the cubic moment gives a three-variable spectral identity on $n$th roots of unity.
Together these identities show that $C$ has no nontrivial $n$th root of unity as a zero.  The
kernel criterion of Hollmann and Xiang then implies $L_c\in\mathscr L_n$.

Combining the two theorems and \cite[Lemma~4.3]{HollmannXiang2004} gives the full structural
classification: for $n>3$, every Kloosterman polynomial is of the form
\[
 L_a\circ L_{1,e},
 \qquad L_a\in\mathscr L_n,
 \qquad e\in\{0,3,6,10\}.
\]

The paper is organized as follows.  \Cref{sec:powersums} evaluates the power sums on $\T$.
\Cref{sec:quadratic} derives the cyclic group-algebra identity supplied by the quadratic
moment.  The cubic moment is developed in \Cref{sec:cubic}, and the central elimination
argument appears in \Cref{sec:elimination}.  The final exclusions proving the necessity in
\Cref{thm:main} are carried out in \Cref{sec:completion}.  The general reduction theorem,
and hence Conjecture~4.5, is proved in \Cref{sec:general-reduction}.

\section{Power sums on the trace-one hyperplane}\label{sec:powersums}

Throughout the proof, indices attached to binary sequences are read in $\Zmod$.  For
\[
 u=\sum_{i=0}^{n-1}u_i2^i,
\]
write
\[
 S(u)=\{i\in\Zmod:u_i=1\}.
\]
The even-weight hypothesis is equivalent to $|S(u)|\equiv0\pmod2$.  We shall therefore classify the even subsets $S\subseteq\Zmod$ for which
\begin{equation}\label{eq:GS}
 G_S(t)=t+\sum_{i\in S}t^{-2^i}
\end{equation}
permutes $\T$.

The cases $n\le3$ are checked directly.  In the necessity argument we henceforth assume $n\ge4$; the small cases will be restored at the end.

Put
\[
 N=2^n-1.
\]
Every element of $\T$ is nonzero, so exponents on $\T$ may be read modulo $N$.

The following power-sum identity is a standard consequence of Newton's identities applied to the roots of the trace polynomial; compare \cite[Theorem~1.75]{LidlNiederreiter1997}.  We include a short Laurent-series proof because the precise residue-class formulation will be used repeatedly.

\begin{lemma}[Power sums on $\T$]\label{lem:powersums}
Assume $n\ge2$.  For every integer $e$,
\[
 \sum_{t\in\T}t^e=
 \begin{cases}
  1,& e\equiv-2^j\pmod N\text{ for some }0\le j<n,\\
  0,& \text{otherwise}.
 \end{cases}
\]
In particular, for $0\le e<N$, the nonzero values occur precisely at
\[
 e=N-1,N-2,N-4,\ldots,N-2^{n-1}.
\]
\end{lemma}

\begin{proof}
The trace-one hyperplane is the zero set of
\[
 A_0(X)=1+X+X^2+X^4+\cdots+X^{2^{n-1}}.
\]
Since $A_0'(X)=1$, one has
\[
 \frac1{A_0(X)}=\frac{A_0'(X)}{A_0(X)}
 =\sum_{t\in\T}\frac1{X+t}.
\]
Expanding at infinity gives
\[
 \frac1{X+t}=\sum_{e\ge0}t^eX^{-e-1},
\]
so the coefficient of $X^{-e-1}$ in $A_0(X)^{-1}$ is $\sum_{t\in\T}t^e$.

On the other hand,
\[
 A_0(X)^2+A_0(X)=X^{2^n}+X,
\]
and therefore
\[
 \frac1{A_0(X)}
 =\frac{A_0(X)+1}{X^{2^n}+X}
 =\frac{\sum_{j=0}^{n-1}X^{2^j}}{X^{2^n}+X}.
\]
As $X^{2^n}+X=X^{2^n}(1+X^{-N})$, we obtain
\[
 \frac1{A_0(X)}
 =\sum_{j=0}^{n-1}\sum_{m\ge0}X^{2^j-2^n-mN}.
\]
The equality $2^j-2^n-mN=-e-1$ is equivalent to
\[
 e=(m+1)N-2^j,
\]
which is precisely the stated congruence.  Reduction of arbitrary $e$ modulo $N$ completes the proof.
\end{proof}

It will be convenient to write
\[
 \Lambda(t^e)=\sum_{t\in\T}t^e.
\]
Thus \Cref{lem:powersums} says that $\Lambda(t^e)=1$ exactly for the residue classes $-2^j$ modulo $N$.

If $G_S$ permutes $\T$, then for every nonnegative integer $m$,
\begin{equation}\label{eq:momentcriterion}
 \sum_{t\in\T}G_S(t)^m=\sum_{t\in\T}t^m.
\end{equation}
We shall use only carefully chosen binary-weight two, three, and four exponents $m$.

\section{The quadratic moment and cyclic orthogonality}\label{sec:quadratic}

Let $s_i\in\F_2$ be the characteristic function of $S$:
\[
 s_i=\begin{cases}1,&i\in S,\\0,&i\notin S.\end{cases}
\]

\begin{proposition}[Cyclic autocorrelation]\label{prop:autocorrelation}
If $G_S$ permutes $\T$, then for every $r\in\Zmod$,
\begin{equation}\label{eq:s-correlation}
 \sum_i s_is_{i+r}=s_{r+1}+s_{1-r}.
\end{equation}
Define
\[
 a_i=s_i+\delta_{i,1},\qquad
 D=S\triangle\{1\},\qquad
 A(Z)=\sum_i a_iZ^i.
\]
Then
\begin{equation}\label{eq:a-correlation}
 \sum_i a_ia_{i+r}=\delta_{r,0},
\end{equation}
and, equivalently,
\begin{equation}\label{eq:unitary-group-ring}
 A(Z)A(Z^{-1})=1
 \quad\text{in}\quad
 \R:=\F_2[Z,Z^{-1}]/(Z^n-1).
\end{equation}
\end{proposition}

\begin{proof}
Use \eqref{eq:momentcriterion} with $m=1+2^r$.  Since
\[
 G_S(t)^{2^r}=t^{2^r}+\sum_{j\in S}t^{-2^{j+r}},
\]
we obtain
\[
 \begin{aligned}
 G_S(t)^{1+2^r}
 ={}&t^{1+2^r}
 +\sum_{j\in S}t^{1-2^{j+r}}
 +\sum_{i\in S}t^{2^r-2^i}\\
 &+\sum_{i,j\in S}t^{-2^i-2^{j+r}}.
 \end{aligned}
\]
The first term cancels the right side of \eqref{eq:momentcriterion}.  By \Cref{lem:powersums}, the second sum contributes only when
\[
 1+2^k\equiv2^{j+r}\pmod N.
\]
The binary uniqueness of a sum of two powers of two gives $k=0$ and $j+r=1$, so the contribution is $s_{1-r}$.  Similarly, the third sum contributes $s_{r+1}$.  Finally,
\[
 -2^i-2^{j+r}\equiv-2^k\pmod N
\]
contributes precisely when $i=j+r$, and hence contributes $\sum_i s_is_{i+r}$.  This proves \eqref{eq:s-correlation}.

Expanding $a_i=s_i+\delta_{i,1}$ in the left side of \eqref{eq:a-correlation} gives the left side of \eqref{eq:s-correlation}, the two correction terms $s_{1-r}$ and $s_{r+1}$, and $\delta_{r,0}$.  Thus \eqref{eq:a-correlation} follows.  Its group-ring formulation is \eqref{eq:unitary-group-ring}, since the coefficient of $Z^r$ in $A(Z)A(Z^{-1})$ is $\sum_i a_ia_{i+r}$.
\end{proof}

The identity \eqref{eq:unitary-group-ring} is the only consequence of the quadratic moment needed below. We now turn to the cubic moments, whose simultaneous form will be combined with \eqref{eq:unitary-group-ring} in the elimination argument.

\section{The cubic moment identity}\label{sec:cubic}

For $p,q\in\Zmod$, define the cubic correlation of the original support by
\[
 T_S(p,q)=\sum_i s_is_{i-p}s_{i-q}.
\]
For $r,s\in\Zmod$, let
\begin{equation}\label{eq:ES-def}
 E_S(r,s)=\Lambda\!\left(
 G_S(t)G_S(t)^{2^r}G_S(t)^{2^s}
 +t^{1+2^r+2^s}
 \right).
\end{equation}
If $G_S$ permutes $\T$, then $E_S(r,s)=0$ for all $r,s$.

\begin{proposition}[Explicit cubic moment]\label{prop:cubic-moment}
For all $r,s\in\Zmod$,
\begin{align}
 E_S(r,s)={}&T_S(r,s-1)+T_S(r-1,s)+T_S(r+1,s+1)\notag\\
 &+s_{-r}s_{-s}+s_rs_{r-s}+s_ss_{s-r}+L_S(r,s),
 \label{eq:ES-explicit}
\end{align}
where
\begin{align}
 L_S(r,s)={}&(\delta_{r,0}+\delta_{r,1})s_{2-s}
 +\delta_{r,-1}s_{1-s}\notag\\
 &+(\delta_{s,0}+\delta_{s,1})s_{2-r}
 +\delta_{s,-1}s_{1-r}\notag\\
 &+\delta_{r,s}s_{r+2}
 +\delta_{s,r+1}s_{r+2}
 +\delta_{r,s+1}s_{s+2}.
 \label{eq:LS}
\end{align}
\end{proposition}

\begin{proof}
All indices in the proof are read in $\Zmod$, and all coefficient computations take place in $\F_2$.  We have
\begin{align*}
 G_S(t)&=t+\sum_i s_it^{-2^i},\\
 G_S(t)^{2^r}&=t^{2^r}+\sum_j s_jt^{-2^{j+r}},\\
 G_S(t)^{2^s}&=t^{2^s}+\sum_{\ell}s_{\ell}t^{-2^{\ell+s}}.
\end{align*}
Thus, in the three factors, the available exponent choices are
\[
\begin{array}{c|c|c}
 \text{factor}&\text{positive exponent}&\text{negative exponent}\\ \hline
 G_S(t)&2^0&-2^i\\
 G_S(t)^{2^r}&2^r&-2^{j+r}\\
 G_S(t)^{2^s}&2^s&-2^{\ell+s}.
\end{array}
\]
The product of the three positive terms is $t^{1+2^r+2^s}$, which is cancelled by the second term in the definition of $E_S(r,s)$.  It remains to examine terms containing one, two, or three negative choices.

We shall repeatedly use two elementary consequences of cyclic binary uniqueness modulo $N=2^n-1$.  First,
\begin{equation}\label{eq:three-power-carry}
 2^a+2^b+2^c\equiv 2^d\pmod N
\end{equation}
holds if and only if, after permuting $a,b,c$, there is a $p\in\Zmod$ such that
\[
 a\equiv p\pmod n,\qquad
 b\equiv p\pmod n,\qquad
 c\equiv p+1\pmod n,\qquad
 d\equiv p+2\pmod n.
\]
This is precisely the cyclic carry
\[
 2^p+2^p+2^{p+1}\equiv2^{p+2}\pmod N.
\]
Second, for fixed $u$,
\begin{equation}\label{eq:two-power-uniqueness}
 2^u+2^v\equiv 2^a+2^b\pmod N
\end{equation}
has a solution $v$ if and only if $a\equiv u\pmod n$ or $b\equiv u\pmod n$.  For completeness, we justify both assertions.  To prove \eqref{eq:three-power-carry}, multiply the congruence by $2^{-d}$ modulo $N$ and choose representatives $0\le \alpha,\beta,\gamma<n$ for $a-d,b-d,c-d$.  Then
\[
 2^\alpha+2^\beta+2^\gamma\equiv1\pmod N.
\]
The left side lies between $3$ and $3\cdot2^{n-1}$, which is strictly smaller than $2N+1$ for $n\ge2$.  Hence it must equal $N+1=2^n$ as an ordinary integer.  After ordering $\alpha\le\beta\le\gamma$, the inequality $\alpha<\beta$ would make
\[
 2^{-\alpha}(2^\alpha+2^\beta+2^\gamma)
 =1+2^{\beta-\alpha}+2^{\gamma-\alpha}
\]
an odd integer greater than one, contrary to the fact that the sum is a power of two.  Thus $\alpha=\beta$.  The equality
\[
 2^{\alpha+1}+2^\gamma=2^n
\]
then forces $\gamma=\alpha+1$ and $n=\alpha+2$.  Translating back gives the pattern $(p,p,p+1)$ with target $p+2$, with all four exponent positions understood modulo $n$.

For \eqref{eq:two-power-uniqueness}, choose all exponents in $\{0,\ldots,n-1\}$.  Each side lies in the interval $[2,2^n]$, so their difference has absolute value less than $N$; the congruence is therefore an equality of ordinary integers.  Ordering the two exponents on each side and comparing their $2$-adic valuations shows that the two multisets of exponents coincide.  In particular, if one exponent on the left is fixed to be $u$, then $a\equiv u\pmod n$ or $b\equiv u\pmod n$.

\smallskip
\noindent\emph{Three negative choices.}
The exponent is
\[
 -2^i-2^{j+r}-2^{\ell+s}.
\]
By \Cref{lem:powersums}, this term contributes precisely when
\[
 2^i+2^{j+r}+2^{\ell+s}\equiv 2^k\pmod N
\]
for some $k$.  By \eqref{eq:three-power-carry}, two of the three displayed exponents must coincide and the third must be their successor.  There are three placements.

If $i=j+r$ and $\ell+s=i+1$, then
\[
 j=i-r,\qquad \ell=i+1-s=i-(s-1),
\]
and the total coefficient is
\[
 \sum_i s_is_{i-r}s_{i-(s-1)}=T_S(r,s-1).
\]
If $i=\ell+s$ and $j+r=i+1$, then
\[
 j=i+1-r=i-(r-1),\qquad \ell=i-s,
\]
and the contribution is
\[
 \sum_i s_is_{i-(r-1)}s_{i-s}=T_S(r-1,s).
\]
Finally, if $j+r=\ell+s$ and $i=j+r+1$, then
\[
 j=i-(r+1),\qquad \ell=i-(s+1),
\]
and the contribution is
\[
 \sum_i s_is_{i-(r+1)}s_{i-(s+1)}=T_S(r+1,s+1).
\]
Hence all-three-negative terms contribute
\begin{equation}\label{eq:three-negative-contribution}
 T_S(r,s-1)+T_S(r-1,s)+T_S(r+1,s+1).
\end{equation}

\smallskip
\noindent\emph{Two negative choices and one positive choice.}
Suppose first that the positive term is selected from $G_S(t)$.  The exponent is
\[
 1-2^{j+r}-2^{\ell+s}.
\]
It is detected by $\Lambda$ exactly when, for some $k$,
\[
 1+2^k\equiv 2^{j+r}+2^{\ell+s}\pmod N.
\]
By \eqref{eq:two-power-uniqueness}, this happens if and only if $j+r=0$ or $\ell+s=0$.  Over $\F_2$, the indicator of a union is
\[
 \ind_{A\lor B}=\ind_A+\ind_B+\ind_A\ind_B.
\]
Therefore this class contributes
\begin{align*}
 &\sum_{j,\ell}s_js_{\ell}
 \bigl(\delta_{j+r,0}+\delta_{\ell+s,0}
       +\delta_{j+r,0}\delta_{\ell+s,0}\bigr)\\
 &\qquad=s_{-r}\sum_{\ell}s_{\ell}
       +s_{-s}\sum_js_j+s_{-r}s_{-s}.
\end{align*}
Since $|S|$ is even, $\sum_i s_i=0$ in $\F_2$, and only $s_{-r}s_{-s}$ remains.

If the positive term comes from $G_S(t)^{2^r}$, the same argument applied to
\[
 2^r-2^i-2^{\ell+s}
\]
shows that the condition is $i=r$ or $\ell+s=r$.  The corresponding coefficient equals
\[
 s_r\sum_{\ell}s_{\ell}+s_{r-s}\sum_i s_i+s_rs_{r-s}
 =s_rs_{r-s}.
\]
Likewise, taking the positive term from $G_S(t)^{2^s}$ gives $s_ss_{s-r}$.  Thus the terms with exactly two negative choices contribute
\begin{equation}\label{eq:two-negative-contribution}
 s_{-r}s_{-s}+s_rs_{r-s}+s_ss_{s-r}.
\end{equation}

\smallskip
\noindent\emph{One negative choice and two positive choices.}
These are the boundary terms encoded by $L_S(r,s)$.  First take the negative term from the third factor.  Its exponent is
\[
 1+2^r-2^{\ell+s},
\]
which is detected precisely when
\begin{equation}\label{eq:one-negative-third}
 1+2^r+2^k\equiv2^{\ell+s}\pmod N
\end{equation}
for some $k$.  Applying \eqref{eq:three-power-carry} to the exponents $0,r,k$ gives exactly three possibilities:
\[
\begin{array}{c|c|c}
 \text{condition on }r&\text{required }k&\ell+s\\ \hline
 r=0&1&2\\
 r=1&0&2\\
 r=-1&-1&1.
\end{array}
\]
After summing over $\ell$, their contribution is
\begin{equation}\label{eq:one-negative-third-contribution}
 (\delta_{r,0}+\delta_{r,1})s_{2-s}
 +\delta_{r,-1}s_{1-s}.
\end{equation}
By symmetry, choosing the negative term from the second factor contributes
\begin{equation}\label{eq:one-negative-second-contribution}
 (\delta_{s,0}+\delta_{s,1})s_{2-r}
 +\delta_{s,-1}s_{1-r}.
\end{equation}

It remains to choose the negative term from the first factor.  The exponent is
\[
 2^r+2^s-2^i,
\]
so the detection condition is
\[
 2^r+2^s+2^k\equiv2^i\pmod N.
\]
Again by \eqref{eq:three-power-carry}, there are precisely three possibilities:
\[
\begin{array}{c|c|c}
 \text{relation between }r,s&\text{required }k&i\\ \hline
 r=s&r+1&r+2\\
 s=r+1&r&r+2\\
 r=s+1&s&s+2.
\end{array}
\]
Consequently this last class contributes
\begin{equation}\label{eq:one-negative-first-contribution}
 \delta_{r,s}s_{r+2}
 +\delta_{s,r+1}s_{r+2}
 +\delta_{r,s+1}s_{s+2}.
\end{equation}
The sum of \eqref{eq:one-negative-third-contribution},
\eqref{eq:one-negative-second-contribution}, and
\eqref{eq:one-negative-first-contribution} is exactly $L_S(r,s)$ in \eqref{eq:LS}.  Combining this with \eqref{eq:three-negative-contribution} and \eqref{eq:two-negative-contribution} proves \eqref{eq:ES-explicit}.
\end{proof}

To exploit this family simultaneously, set
\[
 S(Z)=\sum_i s_iZ^i
\]
and form the two-variable cyclic generating transform
\[
 \widehat E(X,Y)=\sum_{r,s}E_S(r,s)X^rY^s
 \in\F_2[X^{\pm1},Y^{\pm1}]/(X^n-1,Y^n-1).
\]

\begin{corollary}[Two-variable cubic identity]\label{cor:two-variable}
If $G_S$ permutes $\T$, then
\begin{align}
 0={}&(X+Y+X^{-1}Y^{-1})S(XY)S(X^{-1})S(Y^{-1})\notag\\
 &+S(X^{-1})S(Y^{-1})+S(XY)S(Y^{-1})+S(XY)S(X^{-1})\notag\\
 &+\bigl((1+X)Y^2+X^{-1}Y\bigr)S(Y^{-1})\notag\\
 &+\bigl((1+Y)X^2+XY^{-1}\bigr)S(X^{-1})\notag\\
 &+(1+X+Y)X^{-2}Y^{-2}S(XY).
 \label{eq:Ehat}
\end{align}
\end{corollary}

\begin{proof}
Because $G_S$ permutes $\T$, the moment criterion gives $E_S(r,s)=0$ for every $r,s$, and hence $\widehat E(X,Y)=0$.  We compute the transform of the right side of \eqref{eq:ES-explicit} term by term.

First, the basic cubic-correlation transform is
\begin{align}
 \sum_{p,q}T_S(p,q)X^pY^q
 &=\sum_{i,p,q}s_is_{i-p}s_{i-q}X^pY^q\notag\\
 &=\sum_{i,a,b}s_is_as_bX^{i-a}Y^{i-b}\notag\\
 &=S(XY)S(X^{-1})S(Y^{-1}),
 \label{eq:TS-transform}
\end{align}
where we used $a=i-p$ and $b=i-q$.  Put
\[
 F(X,Y)=S(XY)S(X^{-1})S(Y^{-1}).
\]
Changing variables $q=s-1$, $p=r-1$, and $(p,q)=(r+1,s+1)$, respectively, gives
\begin{align*}
 \sum_{r,s}T_S(r,s-1)X^rY^s&=YF(X,Y),\\
 \sum_{r,s}T_S(r-1,s)X^rY^s&=XF(X,Y),\\
 \sum_{r,s}T_S(r+1,s+1)X^rY^s&=X^{-1}Y^{-1}F(X,Y).
\end{align*}
Thus the three cubic-correlation terms transform to
\begin{equation}\label{eq:cubic-transform-total}
 (X+Y+X^{-1}Y^{-1})S(XY)S(X^{-1})S(Y^{-1}).
\end{equation}

For the quadratic terms, direct changes of variables yield
\begin{align}
 \sum_{r,s}s_{-r}s_{-s}X^rY^s
 &=S(X^{-1})S(Y^{-1}),\label{eq:quad-transform-one}\\
 \sum_{r,s}s_rs_{r-s}X^rY^s
 &=S(XY)S(Y^{-1}),\label{eq:quad-transform-two}\\
 \sum_{r,s}s_ss_{s-r}X^rY^s
 &=S(XY)S(X^{-1}).\label{eq:quad-transform-three}
\end{align}
For example, in \eqref{eq:quad-transform-two}, put $a=r$ and $b=r-s$; then $s=a-b$ and
\[
 X^rY^s=(XY)^aY^{-b}.
\]
The third identity is obtained symmetrically.

We next transform the three groups in $L_S(r,s)$.  For the first group,
\begin{align*}
 &\sum_{r,s}\bigl((\delta_{r,0}+\delta_{r,1})s_{2-s}
                  +\delta_{r,-1}s_{1-s}\bigr)X^rY^s\\
 &\quad=(1+X)\sum_s s_{2-s}Y^s
       +X^{-1}\sum_s s_{1-s}Y^s\\
 &\quad=\bigl((1+X)Y^2+X^{-1}Y\bigr)S(Y^{-1}),
\end{align*}
where
\[
 \sum_s s_{2-s}Y^s=Y^2S(Y^{-1}),\qquad
 \sum_s s_{1-s}Y^s=YS(Y^{-1}).
\]
Interchanging $r,X$ with $s,Y$ gives the transform of the second group:
\begin{equation}\label{eq:LS-second-transform}
 \bigl((1+Y)X^2+XY^{-1}\bigr)S(X^{-1}).
\end{equation}

Finally, the three diagonal terms satisfy
\begin{align*}
 \sum_{r,s}\delta_{r,s}s_{r+2}X^rY^s
 &=X^{-2}Y^{-2}S(XY),\\
 \sum_{r,s}\delta_{s,r+1}s_{r+2}X^rY^s
 &=X^{-2}Y^{-1}S(XY),\\
 \sum_{r,s}\delta_{r,s+1}s_{s+2}X^rY^s
 &=X^{-1}Y^{-2}S(XY).
\end{align*}
For instance, in the second line one has $s=r+1$; after setting $i=r+2$,
\[
 X^rY^{r+1}=X^{i-2}Y^{i-1}=X^{-2}Y^{-1}(XY)^i.
\]
Adding the three diagonal transforms gives
\begin{equation}\label{eq:LS-third-transform}
 (1+X+Y)X^{-2}Y^{-2}S(XY).
\end{equation}

Summing \eqref{eq:cubic-transform-total},
\eqref{eq:quad-transform-one}--\eqref{eq:quad-transform-three}, the first transformed group of $L_S$, \eqref{eq:LS-second-transform}, and \eqref{eq:LS-third-transform} yields the right side of \eqref{eq:Ehat}.  Since $\widehat E(X,Y)=0$, the claimed identity follows.
\end{proof}

\section{The group-algebra elimination}\label{sec:elimination}

Recall that $S(Z)=A(Z)+Z$ and $A(Z)A(Z^{-1})=1$.  We now combine \eqref{eq:Ehat} with its symmetry under a change of the three variables whose product is one.

For $pqr=1$, define
\begin{align}
 \E(p,q,r)={}&(p^{-1}+q^{-1}+r^{-1})S(p)S(q)S(r)\notag\\
 &+S(p)S(q)+S(q)S(r)+S(r)S(p)\notag\\
 &+p^{-1}(p^{-1}+q+r)S(p)\notag\\
 &+q^{-1}(q^{-1}+r+p)S(q)\notag\\
 &+r^{-1}(r^{-1}+p+q)S(r).
 \label{eq:calE}
\end{align}
Then \eqref{eq:Ehat} is equivalent to
\begin{equation}\label{eq:E-one}
 \E(X,Y,(XY)^{-1})=0.
\end{equation}
By the symmetry of the three factors in the cubic moment, one also has
\begin{equation}\label{eq:E-two}
 \E(X^{-1},XY,Y^{-1})=0.
\end{equation}

\begin{lemma}[Elimination identity]\label{lem:elimination-identity}
In the two-variable group algebra,
\begin{align}
 &X^3Y^3A(XY)\,\E(X,Y,(XY)^{-1})\notag\\
 &\quad+X^3Y^3A(X)A(Y)\,\E(X^{-1},XY,Y^{-1})\notag\\
 &\qquad=X^2Y^2(X+1)(Y+1)(XY+1)
 \bigl(A(X)A(Y)+A(XY)\bigr).
 \label{eq:master-elimination}
\end{align}
Consequently, if $G_S$ permutes $\T$, then
\begin{equation}\label{eq:factor-H}
 (X+1)(Y+1)(XY+1)\bigl(A(X)A(Y)+A(XY)\bigr)=0.
\end{equation}
\end{lemma}

\begin{proof}
Use $A(U^{-1})=A(U)^{-1}$, which follows from \eqref{eq:unitary-group-ring}, substitute $S(U)=A(U)+U$ into the left side of \eqref{eq:master-elimination}, and clear the monomial denominators.  Equal terms cancel in pairs.  The remaining expression is
\[
 \begin{aligned}
 &(A(X)A(Y)+A(XY))\bigl(X^4Y^4+X^4Y^3+X^3Y^4\\
 &\hspace{42mm}+X^3Y^2+X^2Y^3+X^2Y^2\bigr),
 \end{aligned}
\]
which factors as the right side of \eqref{eq:master-elimination}.  Equations \eqref{eq:E-one} and \eqref{eq:E-two} then give \eqref{eq:factor-H}.
\end{proof}

We now read \eqref{eq:factor-H} coefficientwise.  Put
\[
 H(X,Y)=A(X)A(Y)+A(XY),\qquad
 K(X,Y)=(X+1)(Y+1)H(X,Y).
\]
Then $(XY+1)K(X,Y)=0$.  Write
\[
 K(X,Y)=\sum_{i,j}K_{i,j}X^iY^j.
\]
Comparison of coefficients gives
\begin{equation}\label{eq:diagonal-constant}
 K_{i,j}=K_{i-1,j-1}.
\end{equation}
Thus $K_{i,j}$ is constant along every cyclic diagonal $j-i=d$.

Writing $H(X,Y)=\sum h_{i,j}X^iY^j$, we have
\begin{equation}\label{eq:hcoeff}
 h_{i,j}=\begin{cases}
 a_ia_j,&i\ne j,\\
 0,&i=j.
 \end{cases}
\end{equation}
Since $K=(X+1)(Y+1)H$,
\begin{equation}\label{eq:Kcoeff}
 K_{i,j}=h_{i,j}+h_{i-1,j}+h_{i,j-1}+h_{i-1,j-1}.
\end{equation}
Define the cyclic difference sequence
\begin{equation}\label{eq:c-def}
 c_i=a_i+a_{i-1}.
\end{equation}
If $j-i\notin\{0,1,-1\}$, none of the four entries in \eqref{eq:Kcoeff} lies on the diagonal $i=j$, and therefore
\begin{equation}\label{eq:K-c}
 K_{i,j}=c_ic_j.
\end{equation}
Combining \eqref{eq:diagonal-constant} and \eqref{eq:K-c} gives the rigidity of $c$.

\begin{proposition}[Two-jump rigidity]\label{prop:two-jump}
Assume $n\ge4$ and $G_S$ permutes $\T$.  Then
\[
 A(Z)=Z^k
 \qquad\text{or}\qquad
 A(Z)=1+Z+\cdots+Z^{n-1}+Z^k
\]
for some $k\in\Zmod$.
\end{proposition}

\begin{proof}
For each $d\notin\{0,\pm1\}$, equations \eqref{eq:diagonal-constant} and \eqref{eq:K-c} show that
\begin{equation}\label{eq:cid}
 c_ic_{i+d}\quad\text{is independent of }i.
\end{equation}
This common value must be zero.  Indeed, if it were one for some $d$, then every $c_i$ would be one.  On the adjacent diagonal, a direct use of \eqref{eq:hcoeff} and \eqref{eq:Kcoeff} gives
\begin{equation}\label{eq:adjacent-diagonal}
 K_{i,i+1}=c_ic_{i+1}+a_i.
\end{equation}
The left side is independent of $i$ by \eqref{eq:diagonal-constant}; if $c_i=1$ for every $i$, this forces $a_i$ to be constant, contradicting $c_i=a_i+a_{i-1}=1$.  Hence
\begin{equation}\label{eq:separation-c}
 c_ic_{i+d}=0
 \qquad(i\in\Zmod,\ d\notin\{0,\pm1\}).
\end{equation}

The support of $c$ is therefore contained in two adjacent cyclic positions.  Since
\[
 \sum_i c_i=\sum_i(a_i+a_{i-1})=0,
\]
its weight is even.  Thus either $c=0$ or $\supp(c)=\{t,t+1\}$.

If $c=0$, the sequence $a_i$ is constant.  Equation \eqref{eq:unitary-group-ring} implies $A(1)=1$, so the only possibility is $a_i=1$ for all $i$ and $n$ odd.  But with
\[
 N_n(Z)=1+Z+\cdots+Z^{n-1}
\]
one has $N_n(Z)N_n(Z^{-1})=N_n(Z)\ne1$ for $n>1$, contradicting \eqref{eq:unitary-group-ring}.

Hence $\supp(c)=\{t,t+1\}$.  A cyclic binary sequence with exactly two adjacent jumps has either a single coefficient equal to one or a single coefficient equal to zero.  This gives respectively
\[
 A(Z)=Z^k
 \quad\text{or}\quad
 A(Z)=N_n(Z)+Z^k.
\]
\end{proof}

The first alternative is the desired monomial conclusion.  It remains to rule out the complementary monomial.

\section{Completion of the classification}\label{sec:completion}

\subsection{Excluding the complementary monomial}

Suppose
\begin{equation}\label{eq:complementary}
 A(Z)=N_n(Z)+Z^k.
\end{equation}
Since $A(1)=1$, this can occur only when $n$ is even.  Moreover
\[
 S(Z)=A(Z)+Z=N_n(Z)+Z^k+Z.
\]
For a residue $i\in\Zmod$, write
\[
 [i]=\begin{cases}1,&k\equiv i\pmod n,\\0,&\text{otherwise}.
 \end{cases}
\]
Then
\begin{equation}\label{eq:s-complement}
 s_i=1+[i]+\delta_{i,1}.
\end{equation}
We now give the calculation of the three cubic moments used to exclude this case.  Put
\[
 v_i=\delta_{i,1}+\delta_{i,k},
 \qquad s_i=1+v_i.
\]
Because $n$ is even and the auxiliary word $v$ has either zero or two nonzero coordinates,
\[
 \sum_i1=0,
 \qquad
 \sum_i v_i=0
 \qquad\text{in }\F_2.
\]
Expanding the cubic correlation gives
\begin{equation}\label{eq:T-complement-expansion}
 T_S(p,q)=R(p)+R(q)+R(q-p)+U(p,q),
\end{equation}
where
\[
 R(d)=\sum_i v_i v_{i-d},
 \qquad
 U(p,q)=\sum_i v_i v_{i-p}v_{i-q}.
\]
Since $v_i$ is supported at the two cyclic positions $1$ and $k$, one has
\begin{equation}\label{eq:R-complement}
 R(d)=[d+1]+[1-d].
\end{equation}
The triple-intersection term $U(p,q)$ is obtained by requiring all three positions
$i,i-p,i-q$ to lie in $\{1,k\}$.  Applying this observation to the values needed below yields
\begin{align}
 T_S(1,1)&=[0]+[2],
 &T_S(0,2)&=[-1]+[3],\notag\\
 T_S(2,3)&=[-2]+[-1]+[0]+[2]+[3]+[4],
 &T_S(1,2)&=[-1]+[3],\notag\\
 T_S(0,3)&=[-2]+[4],
 &T_S(2,4)&=[-3]+[5],\notag\\
 T_S(1,4)&=[-3]+[-2]+[0]+[2]+[4]+[5],\notag\\
 T_S(3,5)&=[-4]+[-2]+[-1]+[3]+[4]+[6].
 \label{eq:T-complement-values}
\end{align}
All identities in this subsection are in $\F_2$, and every index is read modulo $n$.
For example, \eqref{eq:T-complement-expansion} gives
\[
 T_S(1,1)=R(0)+R(1)+R(1)+U(1,1)=U(1,1)=[0]+[2].
\]
The other entries in \eqref{eq:T-complement-values} follow in exactly the same way from
\eqref{eq:R-complement} and the two-point support of $v$.

We first take $(r,s)=(1,2)$ in \eqref{eq:ES-explicit}.  The boundary term is
\[
 L_S(1,2)=s_0+s_3,
\]
and hence
\begin{align*}
 E_S(1,2)={}&T_S(1,1)+T_S(0,2)+T_S(2,3)\\
 &+s_{-1}s_{-2}+s_1s_{-1}+s_2s_1+s_0+s_3.
\end{align*}
By \eqref{eq:T-complement-values}, the sum of the three cubic correlations is
\[
 T_S(1,1)+T_S(0,2)+T_S(2,3)=[-2]+[4].
\]
Since $n\ge10$, the cyclic positions occurring in each product below are distinct unless they
are explicitly equal.  Thus
\[
 s_{-1}s_{-2}=1+[-1]+[-2],
 \qquad
 s_1s_{-1}=[1],
 \qquad
 s_2s_1=[1],
\]
and
\[
 s_0+s_3=(1+[0])+(1+[3])=[0]+[3].
\]
Consequently
\begin{equation}\label{eq:comp1}
 \begin{aligned}
 E_S(1,2)
 &=([-2]+[4])+(1+[-1]+[-2])+([0]+[3])\\
 &=1+[0]+[3]+[4]+[-1].
 \end{aligned}
\end{equation}

Next take $(r,s)=(1,3)$.  Here $L_S(1,3)=s_{-1}$, and therefore
\begin{align*}
 E_S(1,3)={}&T_S(1,2)+T_S(0,3)+T_S(2,4)\\
 &+s_{-1}s_{-3}+s_1s_{-2}+s_3s_2+s_{-1}.
\end{align*}
The cubic-correlation contribution is
\[
 T_S(1,2)+T_S(0,3)+T_S(2,4)
 =[-3]+[-2]+[-1]+[3]+[4]+[5].
\]
Moreover,
\[
 s_{-1}s_{-3}=1+[-1]+[-3],
 \qquad
 s_1s_{-2}=[1],
 \qquad
 s_3s_2=1+[3]+[2],
\]
while $s_{-1}=1+[-1]$.  The last four terms consequently sum to
\[
 1+[-3]+[1]+[2]+[3].
\]
It follows that
\begin{equation}\label{eq:comp2}
 \begin{aligned}
 E_S(1,3)
 &=[-3]+[-2]+[-1]+[3]+[4]+[5]\\
 &\quad+1+[-3]+[1]+[2]+[3]\\
 &=1+[1]+[2]+[4]+[5]+[-2]+[-1].
 \end{aligned}
\end{equation}

Finally take $(r,s)=(2,4)$.  For $n\ge10$, none of the Kronecker-delta conditions in
\eqref{eq:LS} occurs, so
\[
 L_S(2,4)=0.
\]
Thus
\begin{align*}
 E_S(2,4)={}&T_S(2,3)+T_S(1,4)+T_S(3,5)\\
 &+s_{-2}s_{-4}+s_2s_{-2}+s_4s_2.
\end{align*}
The first line contributes
\[
 T_S(2,3)+T_S(1,4)+T_S(3,5)
 =[-4]+[-3]+[-2]+[4]+[5]+[6].
\]
For the quadratic terms,
\[
 s_{-2}s_{-4}=1+[-2]+[-4],
 \qquad
 s_2s_{-2}=1+[2]+[-2],
 \qquad
 s_4s_2=1+[4]+[2].
\]
Their sum is $1+[-4]+[4]$, and hence
\begin{equation}\label{eq:comp3}
 \begin{aligned}
 E_S(2,4)
 &=([-4]+[-3]+[-2]+[4]+[5]+[6])\\
 &\quad +(1+[-4]+[4])\\
 &=1+[5]+[6]+[-3]+[-2].
 \end{aligned}
\end{equation}

If all three cubic moments vanish, \eqref{eq:comp1} first forces
\[
 k\in\{0,3,4,-1\}.
\]
Intersecting with the vanishing set of \eqref{eq:comp2} leaves only
\[
 k\in\{4,-1\}.
\]
For either of these residues none of the four indicators in \eqref{eq:comp3} is one, so
$E_S(2,4)=1$, a contradiction.  Thus the complementary monomial cannot occur for even
$n\ge10$.

The remaining even values $n=4,6,8$ involve cyclic coincidences among the displayed residues and are read directly from \eqref{eq:ES-explicit}.  The resulting checks are summarized in \Cref{tab:complement-small}.

\begin{table}[ht]
\centering
\caption{Finite checks for the complementary monomial.}\label{tab:complement-small}
\begin{tabular}{ccl}
\toprule
$n$ & moment(s) & obstruction\\
\midrule
$4$ & $E_S(1,2)$ & equals $1$ for every $k$\\
$6$ & $E_S(1,2),E_S(1,3)$ & their vanishing sets are disjoint\\
$8$ & $E_S(2,4)$ & equals $1$ for every $k$\\
\bottomrule
\end{tabular}
\end{table}

We have proved the following.

\begin{corollary}\label{cor:monomial-A}
If $n\ge4$ and $G_S$ permutes $\T$, then
\[
 A(Z)=Z^k
\]
for some $k\in\Zmod$.  Consequently
\begin{equation}\label{eq:S-two-term}
 S=\{1,k\},
\end{equation}
with the convention that $S=\varnothing$ when $k=1$.
\end{corollary}

\subsection{A cubic restriction on the remaining exponent}

By \eqref{eq:S-two-term}, the map has the form
\begin{equation}\label{eq:Gk}
 G_k(t)=t+t^{-2}+t^{-2^k}.
\end{equation}
The values $k=1,0,2,3$ correspond respectively to
\[
 S=\varnothing,\quad \{0,1\},\quad \{1,2\},\quad \{1,3\},
\]
and hence to $u=0,3,6,10$.  We now exclude $k\ge4$.

\begin{lemma}[Cubic restriction]\label{lem:cubic-restriction}
Assume $k\ge4$ and $G_k$ permutes $\T$.  Then
\begin{equation}\label{eq:2k-congruence}
 2k\equiv3\pmod n.
\end{equation}
In particular, after choosing $0\le k<n$, one has
\[
 n\text{ odd},\qquad k=\frac{n+3}{2}.
\]
\end{lemma}

\begin{proof}
Since $0\le k<n$ and $k\ge4$, we have $n\ge5$.  Take the moment exponent
\[
 m=2^{k-1}+3=2^{k-1}+2+1.
\]
In characteristic two,
\[
 G_k(t)^m=G_k(t)\,G_k(t)^2\,G_k(t)^{2^{k-1}}.
\]
For the factor $G_k(t)^{2^r}$, write its three available exponents as
\[
 P_r=2^r,\qquad B_r=-2^{r+1},\qquad C_r=-2^{r+k}.
\]
Thus the three factors, corresponding to $r=0,1,k-1$, offer the choices
\[
 \{P_0,B_0,C_0\},\qquad
 \{P_1,B_1,C_1\},\qquad
 \{P_{k-1},B_{k-1},C_{k-1}\}.
\]
There are $27$ resulting monomials.  The all-positive choice has exponent
\[
 P_0+P_1+P_{k-1}=1+2+2^{k-1}=m,
\]
so it cancels the right-hand side of the moment identity
\[
 \sum_{t\in\T}G_k(t)^m=\sum_{t\in\T}t^m.
\]
We therefore count, in $\F_2$, the contributions of the remaining $26$ choices.  By
\Cref{lem:powersums}, a choice contributes precisely when its exponent sum is congruent to
$-2^j$ modulo $N=2^n-1$ for some $j\in\Zmod$.

Three choices contribute for every $n$ and $k$:
\begin{align*}
 B_0+P_1+B_{k-1}
   &=(-2)+2-2^k=-2^k,\\
 B_0+P_1+C_{k-1}
   &=(-2)+2-2^{2k-1}=-2^{2k-1},\\
 C_0+C_1+B_{k-1}
   &=-2^k-2^{k+1}-2^k=-2^{k+2}.
\end{align*}
Hence these universal terms contribute $3\equiv1\pmod 2$.

For completeness, we now record the entire list of choices that can contribute in addition to
these three.  Juxtaposition indicates one choice from each of the three factors.  Repeated use of
the cyclic three-power carry \eqref{eq:three-power-carry} and the two-power uniqueness
\eqref{eq:two-power-uniqueness} gives the following table; the other $17$ choices never satisfy
the power-sum criterion.
\begin{table}[ht]
\centering
\renewcommand{\arraystretch}{1.18}
\begin{tabular}{@{}lll@{}}
\toprule
Choice & Exponent sum & Necessary and sufficient condition \\
\midrule
$P_0P_1C_{k-1}$ & $3-2^{2k-1}$
  & $2k\equiv3\pmod n$ \\
$P_0B_1C_{k-1}$ & $-3-2^{2k-1}$
  & $2k\equiv1\pmod n$ \\
$P_0C_1B_{k-1}$ & $1-2^{k+1}-2^k$
  & $k\equiv-1\pmod n$ \\
$P_0C_1C_{k-1}$ & $1-2^{k+1}-2^{2k-1}$
  & $k\equiv-1\pmod n$ or $2k\equiv1\pmod n$ \\
$B_0B_1C_{k-1}$ & $-2-4-2^{2k-1}$
  & $2k\equiv2\pmod n$ \\
$C_0P_1C_{k-1}$ & $2-2^k-2^{2k-1}$
  & $2k\equiv2\pmod n$ \\
\bottomrule
\end{tabular}
\caption{All possible non-universal corrections in the cubic moment.}
\label{tab:cubic-corrections}
\end{table}

We include the carry checks behind the table.  For the first row, contribution is equivalent to
\[
 3+2^j\equiv2^{2k-1}\pmod N.
\]
The left side is $1+2+2^j$.  By \eqref{eq:three-power-carry}, it can be a single power of two
only when $j\equiv0\pmod n$, in which case $1+2+1=4$.  Hence
\[
 2^{2k-1}\equiv2^2\pmod N,
 \qquad\text{or equivalently}\qquad
 2k\equiv3\pmod n.
\]
Similarly, the second and fifth rows reduce respectively to
\[
 1+2+2^{2k-1}\equiv2^j\pmod N,
 \qquad
 2+4+2^{2k-1}\equiv2^j\pmod N.
\]
The three-power carry forces, respectively,
\[
 2k-1\equiv0\pmod n,
 \qquad
 2k-1\equiv1\pmod n,
\]
which give $2k\equiv1\pmod n$ and $2k\equiv2\pmod n$.

The remaining rows follow from two-power uniqueness.  The third row contributes exactly when
\[
 1+2^j\equiv2^k+2^{k+1}\pmod N.
\]
Thus the exponent $0$ on the left must satisfy $0\equiv k\pmod n$ or $0\equiv k+1\pmod n$.  Since $0\le k<n$ and
$k\ge4$, the first alternative is impossible, and hence $k+1\equiv0\pmod n$, that is,
$k\equiv-1\pmod n$.  For the fourth row one obtains
\[
 1+2^j\equiv2^{k+1}+2^{2k-1}\pmod N,
\]
so $0\equiv k+1\pmod n$ or $0\equiv2k-1\pmod n$.  This gives
\[
 k\equiv-1\pmod n
 \qquad\text{or}\qquad
 2k\equiv1\pmod n.
\]
Finally, the sixth row gives
\[
 2+2^j\equiv2^k+2^{2k-1}\pmod N.
\]
The fixed exponent $1$ must therefore satisfy $1\equiv k\pmod n$ or $1\equiv2k-1\pmod n$.  Again $k\equiv1\pmod n$ is impossible,
so $2k-1\equiv1\pmod n$, namely $2k\equiv2\pmod n$.

It remains to count these corrections modulo two.  Under the condition $2k\equiv1\pmod n$, the two
terms
\[
 P_0B_1C_{k-1},\qquad P_0C_1C_{k-1}
\]
occur together.  Under $k\equiv-1\pmod n$, the terms
\[
 P_0C_1B_{k-1},\qquad P_0C_1C_{k-1}
\]
occur together.  Under $2k\equiv2\pmod n$, the terms
\[
 B_0B_1C_{k-1},\qquad C_0P_1C_{k-1}
\]
occur together.  These three side conditions are pairwise incompatible for $n\ge5$: an overlap
would force $n$ to divide $1$, $3$, or $4$.  Consequently all corrections in these three classes
occur in disjoint pairs and make zero total contribution in $\F_2$.

The only correction not forced to occur in a pair is $P_0P_1C_{k-1}$, whose indicator is
$\ind_{\{2k\equiv3\pmod n\}}$.  After cancellation of the all-positive term, the parity of the
entire cubic moment is therefore
\[
 1+\ind_{\{2k\equiv3\pmod n\}}.
\]
The moment identity requires this quantity to vanish.  Hence
\[
 2k\equiv3\pmod n,
\]
which proves \eqref{eq:2k-congruence}.

Finally, $0\le k<n$ and $k\ge4$ give
\[
 0<2k-3<2n.
\]
Since $n$ divides $2k-3$, it follows that $2k-3=n$.  Thus $n$ is odd and
\[
 k=\frac{n+3}{2}.
\]
\end{proof}

\subsection{The quartic moment exclusion}

Assume now that
\begin{equation}\label{eq:k-special}
 n\text{ is odd},\qquad k=\frac{n+3}{2}.
\end{equation}
Take
\[
 m=2^{k-1}+7=2^{k-1}+4+2+1
\]
and put $K=2^{k-1}$.  The four factors correspond to $r=0,1,2,k-1$.  Since $2k\equiv3\pmod n$, the exponent $-2^{2k-1}$ is congruent to $-4$ modulo $N$.  The choices in the four factors are therefore represented by
\begin{equation}\label{eq:quartic-table}
\begin{array}{c|ccc}
 r & P_r & B_r & C_r\\
\hline
0 & 1 & -2 & -2K\\
1 & 2 & -4 & -4K\\
2 & 4 & -8 & -8K\\
k-1 & K & -2K & -4.
\end{array}
\end{equation}

\begin{lemma}[Quartic exclusion]\label{lem:quartic-exclusion}
No $k\ge4$ satisfying \eqref{eq:k-special} can make $G_k$ a permutation of $\T$.
\end{lemma}

\begin{proof}
We first treat the stable range $n\ge11$.  Since $n=2k-3$, this is equivalent to
$k\ge7$, and hence
\[
 K=2^{k-1}\ge64.
\]
In terms of the decomposition into a multiple of $K$ and a bounded remainder, the four rows
of \eqref{eq:quartic-table} are
\begin{equation}\label{eq:quartic-HL-table}
\begin{array}{c|ccc}
 r & P_r & B_r & C_r\\
\hline
0   &(0)K+1 &(0)K-2 &(-2)K\\
1   &(0)K+2 &(0)K-4 &(-4)K\\
2   &(0)K+4 &(0)K-8 &(-8)K\\
k-1 &(1)K   &(-2)K  &(0)K-4.
\end{array}
\end{equation}
Thus every exponent occurring in the expansion of
\[
 G_k(t)^{2^{k-1}+7}
   =G_k(t)G_k(t)^2G_k(t)^4G_k(t)^{2^{k-1}}
\]
has a representation
\begin{equation}\label{eq:E-HK-L}
 E=HK+L,
\end{equation}
where $H,L\in\mathbb Z$ and
\begin{equation}\label{eq:L-bound}
 |L|\le18<\frac K2.
\end{equation}
The representation arising from \eqref{eq:quartic-HL-table} is unique: if
$HK+L=H'K+L'$ with $|L|,|L'|\le18$, then
\[
 |(H-H')K|=|L'-L|\le36<K,
\]
so $H=H'$ and $L=L'$.

We first justify carefully why the power-sum congruence can be treated as an ordinary integer
equality in this range.  Put
\[
 N=2^n-1=2^{2k-3}-1=2^{k-2}K-1.
\]
Every exponent in \eqref{eq:quartic-HL-table} satisfies
\begin{equation}\label{eq:E-range}
 -16K\le E\le K+7.
\end{equation}
Since $k\ge7$, we have $N\ge32K-1$, and therefore $-N<E<N$.  Choose
$0\le j<n$.  If
\[
 E\equiv-2^j\pmod N,
\]
then, because both $E$ and $-2^j$ lie in $(-N,N)$, the only possible integer equalities are
\[
 E=-2^j
 \qquad\text{or}\qquad
 E=N-2^j.
\]
The second alternative is impossible.  Indeed,
\[
 N-2^j\ge N-2^{n-1}=2^{n-1}-1
      =2^{k-3}K-1>K+7,
\]
contradicting \eqref{eq:E-range}.  Consequently
\begin{equation}\label{eq:congruence-to-equality}
 E\equiv-2^j\pmod N
 \quad\Longleftrightarrow\quad
 E=-2^j
\end{equation}
for every exponent $E$ appearing in the quartic expansion.

We now use \eqref{eq:E-HK-L}.  Suppose first that $j\ge k-1$.  Then
\[
 -2^j=-2^{j-k+1}K
\]
is a multiple of $K$.  Reducing $HK+L=-2^j$ modulo $K$ gives
$L\equiv0\pmod K$.  Since $|L|<K$, it follows that
\begin{equation}\label{eq:large-j-HL}
 L=0,
 \qquad
 H=-2^{j-k+1}.
\end{equation}
Thus $H$ must be a negative power of two.

Suppose instead that $j<k-1$.  Then $2^j\le K/2$.  If $H\ge1$, then
$HK+L\ge K-18>0$, whereas $-2^j<0$.  If $H\le-2$, then
$HK+L\le-2K+18<-K/2\le-2^j$.  The only remaining nonzero possibility is
$H=-1$, but then
\[
 L=K-2^j\ge\frac K2>18,
\]
contrary to \eqref{eq:L-bound}.  Hence
\begin{equation}\label{eq:small-j-HL}
 H=0,
 \qquad
 L=-2^j.
\end{equation}
This explains precisely the two alternatives used below.

We next enumerate the possible remainders.  Selecting one entry from each row of
\eqref{eq:quartic-HL-table}, subject to $H=0$, gives
\begin{equation}\label{eq:H-zero-L-values}
 L\in\{-18,-15,-12,-9,-6,-3,0,3\}.
\end{equation}
None of these values belongs to
\[
 \{-1,-2,-4,-8,-16\},
\]
the list of negative powers of two with absolute value at most $18$.  Therefore the case
$j<k-1$ produces no contribution.

Similarly, imposing $L=0$ gives
\begin{equation}\label{eq:L-zero-H-values}
 H\in\{-16,-13,-10,-7,-6,-4,-1,0\}.
\end{equation}
The negative powers of two in this list are exactly $-1,-4,-16$.  They arise from precisely
the following three choices:
\begin{align}
 C_0+B_1+P_2+P_{k-1}
   &=-2K-4+4+K=-K=-2^{k-1},\label{eq:quartic-hit-1}\\
 C_0+B_1+P_2+B_{k-1}
   &=-2K-4+4-2K=-4K=-2^{k+1},\label{eq:quartic-hit-2}\\
 C_0+C_1+C_2+B_{k-1}
   &=-2K-4K-8K-2K=-16K=-2^{k+3}.
   \label{eq:quartic-hit-3}
\end{align}
For $k\ge7$, the exponents $k-1,k+1,k+3$ all lie between $0$ and $n-1=2k-4$.
Thus these three and only these three monomials are detected by \Cref{lem:powersums}.  Hence
\begin{equation}\label{eq:quartic-left-stable}
 \sum_{t\in\T}G_k(t)^{2^{k-1}+7}=1
 \qquad\text{in }\F_2.
\end{equation}

It remains to evaluate the corresponding moment of the identity map.  We claim that, for every
$k\ge5$ with $n=2k-3$,
\begin{equation}\label{eq:quartic-rhs-nonhit}
 2^{k-1}+7\not\equiv-2^j\pmod{2^n-1}
 \qquad(0\le j<n).
\end{equation}
Indeed, $0<2^{k-1}+7<N$, so such a congruence would imply
\[
 2^{k-1}+7=N-2^j=2^{2k-3}-1-2^j,
\]
or equivalently
\begin{equation}\label{eq:quartic-binary-impossible}
 2^{2k-3}=2^j+2^{k-1}+2^3.
\end{equation}
The right-hand side cannot be a single power of two.  If
$j\notin\{3,k-1\}$, its binary expansion has three distinct nonzero positions.  If $j=3$, then
it becomes
\[
 2^4+2^{k-1};
\]
for $k>5$ the two positions are distinct, while for $k=5$ the value is $2^5$, not
$2^{2k-3}=2^7$.  If $j=k-1$, then it becomes
\[
 2^k+2^3,
\]
which again has two distinct nonzero binary positions because $k\ge5$.  This proves
\eqref{eq:quartic-rhs-nonhit}, and therefore
\begin{equation}\label{eq:quartic-right-stable}
 \sum_{t\in\T}t^{2^{k-1}+7}=0
 \qquad(k\ge5).
\end{equation}
Equations \eqref{eq:quartic-left-stable} and \eqref{eq:quartic-right-stable} contradict the
moment criterion, completing the proof for $n\ge11$.

It remains to inspect the residual pairs
\[
 (n,k)=(5,4),(7,5),(9,6).
\]
The same $81$-term expansion gives the following explicit lists.  When $(n,k)=(9,6)$, the
only detected choices are the three choices in
\eqref{eq:quartic-hit-1}--\eqref{eq:quartic-hit-3}; hence the left-hand moment is $1$, while
the right-hand moment is zero by \eqref{eq:quartic-rhs-nonhit}.

When $(n,k)=(7,5)$, exactly the following five choices are detected:
\begin{equation}\label{eq:n7-hits}
 \begin{gathered}
 P_0P_1C_2C_{k-1},\qquad
 P_0B_1C_2C_{k-1},\\
 C_0B_1P_2P_{k-1},\qquad
 C_0B_1P_2B_{k-1},\qquad
 C_0C_1C_2B_{k-1}.
 \end{gathered}
\end{equation}
Thus the left-hand moment is $1$, whereas the right-hand moment is again zero.

Finally, let $(n,k)=(5,4)$.  The eight detected choices are
\begin{equation}\label{eq:n5-hits}
 \begin{gathered}
 P_0P_1P_2P_{k-1},\qquad
 C_0P_1P_2P_{k-1},\\
 C_0P_1C_2P_{k-1},\qquad
 C_0P_1C_2B_{k-1},\\
 C_0B_1P_2P_{k-1},\qquad
 C_0B_1P_2B_{k-1},\\
 C_0B_1B_2C_{k-1},\qquad
 C_0C_1C_2B_{k-1}.
 \end{gathered}
\end{equation}
Consequently the left-hand moment is $0$.  On the other hand,
\[
 2^{k-1}+7=15\equiv-16=-2^4\pmod{31},
\]
so the right-hand moment is $1$.  Each of the three residual pairs therefore contradicts
\eqref{eq:momentcriterion}.  No $k\ge4$ satisfying \eqref{eq:k-special} can survive.
\end{proof}

Combining \Cref{cor:monomial-A} with \Cref{lem:cubic-restriction,lem:quartic-exclusion} proves that, for $n\ge4$, the only possible supports are
\[
 \varnothing,\qquad \{0,1\},\qquad \{1,2\},\qquad \{1,3\}.
\]
Equivalently,
\[
 u\in\{0,3,6,10\}.
\]
Together with the direct cases $n\le3$, this proves the necessity in \Cref{thm:main}.  The sufficiency is the trivial case $u=0$ together with \cite[Theorem~4.1]{HollmannXiang2004}.  Hence \Cref{thm:main}, and therefore Conjecture~4.4, is proved.

\section{The general Hollmann--Xiang reduction}\label{sec:general-reduction}

We now prove \Cref{thm:general-reduction}, which is
\cite[Conjecture~4.5]{HollmannXiang2004}.  Throughout this section let $n>3$, put
\[
 F=\F_{2^n},\qquad
 \mathcal T_e=\{x\in F:\Tr(x)=e\}\quad(e\in\F_2),
\]
and write
\[
 C(Z)=\sum_{i=0}^{n-1}c_iZ^i,
 \qquad
 D(Z)=\sum_{i=0}^{n-1}d_iZ^i.
\]
The corresponding linearized maps are
\[
 L_C(x)=\sum_i c_ix^{2^i},
 \qquad
 L_D(x)=\sum_i d_ix^{2^i}.
\]
Assume that $D(1)=0$ and that
\begin{equation}\label{eq:general-f}
 f(x)=L_C(x)+L_D(x^{-1})
\end{equation}
is injective on $\mathcal T_1$.  If
\[
 \varepsilon=C(1)\in\F_2,
\]
then
\[
 \Tr(f(x))=C(1)\Tr(x)+D(1)\Tr(x^{-1})=\varepsilon
 \qquad(x\in\mathcal T_1).
\]
Since both trace fibers have cardinality $2^{n-1}$, the map in
\eqref{eq:general-f} is therefore a bijection
\begin{equation}\label{eq:f-bijection}
 f:\mathcal T_1\longrightarrow\mathcal T_\varepsilon.
\end{equation}

We first record the linearized kernel criterion used at the end of the argument.  It is
\cite[Lemma~4.3 and Remark~1]{HollmannXiang2004} in polynomial form.

\begin{lemma}[Trace-fiber kernel criterion]\label{lem:trace-fiber-kernel}
Let
\[
 H(Z)=\sum_{i=0}^{n-1}h_iZ^i\in\F_2[Z],
 \qquad
 L_H(x)=\sum_{i=0}^{n-1}h_ix^{2^i}.
\]
Then $L_H$ is injective on $\mathcal T_0$ if and only if it is injective on
$\mathcal T_1$.  These equivalent conditions hold precisely in one of the following cases:
\begin{enumerate}[label=\textup{(\alph*)}]
 \item $H(1)=1$ and
 \[
  \gcd(H(Z),Z^n-1)=1;
 \]
 \item $H(1)=0$, the integer $n$ is odd, and
 \[
  \gcd(H(Z),Z^n-1)=Z+1.
 \]
\end{enumerate}
\end{lemma}

\begin{proof}
Let $K=\ker L_H$.  For $x,y$ in the same trace fiber, one has
$L_H(x)=L_H(y)$ if and only if $x+y\in K\cap\mathcal T_0$.  Hence injectivity on either
trace fiber is equivalent to
\[
 K\cap\mathcal T_0=\{0\}.
\]
If $K=\{0\}$, then $L_H$ is a permutation of $F$ and necessarily $H(1)=1$.  If $K$ is
nontrivial and has dimension at least two, then the hyperplane $\mathcal T_0$ meets $K$
nontrivially.  Thus the only other possibility is
\[
 K=\{0,b\}
\]
with $\Tr(b)=1$.  Since
\[
 \Tr(L_H(x))=H(1)\Tr(x),
\]
this case forces $H(1)=0$, whence $L_H(1)=0$ and therefore $b=1$.  It also forces
$\Tr(1)=1$, so $n$ is odd.  Conversely, each of these two kernel configurations gives
injectivity on both trace fibers.

Finally, the standard normal-basis criterion for binary linearized polynomials identifies the
$\F_2$-dimension of $\ker L_H$ with the degree of $\gcd(H(Z),Z^n-1)$.  Thus $K=\{0\}$ is
equivalent to the gcd being $1$, while $K=\{0,1\}$ is equivalent, in the odd-degree case, to
the gcd being $Z+1$.
\end{proof}

\subsection{The parity obstruction}

\begin{lemma}\label{lem:even-degree-parity}
If $n$ is even, then $C(1)=1$.
\end{lemma}

\begin{proof}
Suppose instead that $C(1)=0$.  Then \eqref{eq:f-bijection} is a bijection from
$\mathcal T_1$ onto $\mathcal T_0$.  Hence there is a unique $x_0\in\mathcal T_1$ such that
$f(x_0)=0$.  Since all coefficients lie in $\F_2$,
\[
 f(x^2)=f(x)^2.
\]
Thus $f(x_0^2)=0$.  The trace is Frobenius invariant, so $x_0^2\in\mathcal T_1$, and
uniqueness gives $x_0^2=x_0$.  Therefore $x_0\in\F_2$.  But if $n$ is even, then both
$0$ and $1$ have trace zero, contradicting $x_0\in\mathcal T_1$.
\end{proof}

In particular, the theorem is already proved for $n=4$: in that case
\[
 Z^4-1=(Z+1)^4,
\]
and \Cref{lem:even-degree-parity} gives $C(1)=1$, so
$\gcd(C(Z),Z^4-1)=1$.  We may therefore assume from now on that
\begin{equation}\label{eq:general-n-ge5}
 n\ge5.
\end{equation}

\subsection{The relative quadratic moment}

Put
\[
 N=2^n-1,
 \qquad
 \mathcal R_n=\F_2[Z,Z^{-1}]/(Z^n-1),
\]
and define
\begin{equation}\label{eq:def-B-general}
 B(Z)=D(Z)+ZC(Z)\in\mathcal R_n.
\end{equation}
If $B(Z)=\sum_i b_iZ^i$, then
\[
 b_i=d_i+c_{i-1},
\]
with all subscripts read modulo $n$.

\begin{proposition}[Relative quadratic identity]\label{prop:relative-quadratic}
In $\mathcal R_n$ one has
\begin{equation}\label{eq:relative-quadratic}
 B(Z)B(Z^{-1})=C(Z)C(Z^{-1}).
\end{equation}
Consequently, for every $z\in\overline{\F}_2$ satisfying $z^n=1$,
\begin{equation}\label{eq:relative-quadratic-eval}
 B(z)B(z^{-1})=C(z)C(z^{-1}).
\end{equation}
\end{proposition}

\begin{proof}
Fix $r\in\mathbb Z/n\mathbb Z$.  Since $f$ is a bijection as in
\eqref{eq:f-bijection}, its $(1+2^r)$th power sum satisfies
\begin{equation}\label{eq:general-second-moment}
 \sum_{x\in\mathcal T_1}f(x)^{1+2^r}
 =\sum_{y\in\mathcal T_\varepsilon}y^{1+2^r}.
\end{equation}
The right-hand side is zero.  For $\mathcal T_1$, this follows from
\Cref{lem:powersums}, because $1+2^r$ is not congruent to $-2^j$ modulo
$N=2^n-1$ when $n>3$.  For $\mathcal T_0$, use
$F=\mathcal T_0\sqcup\mathcal T_1$: the sum of $y^{1+2^r}$ over $F$ is zero,
and so the sums over the two trace fibers are equal and hence both vanish.

We now expand the left-hand side.  Write
\[
 f(x)=\sum_i c_ix^{2^i}+\sum_i d_ix^{-2^i},
 \qquad
 f(x)^{2^r}=\sum_j c_jx^{2^{j+r}}+\sum_jd_jx^{-2^{j+r}}.
\]
All indices are read modulo $n$.  A term contributes to the sum over
$\mathcal T_1$ precisely when its exponent is congruent to $-2^k$ modulo
$N$ for some $k$.

First choose a negative term from both factors.  The exponent is
\[
 -2^i-2^{j+r}.
\]
It is detected precisely when
\[
 2^i+2^{j+r}\equiv2^k\pmod N.
\]
The only binary carry producing a single power of two is obtained when
$i\equiv j+r\pmod n$, in which case $k\equiv i+1\pmod n$.  Thus the
contribution is
\[
 \sum_i d_id_{i-r}=\sum_i d_id_{i+r}.
\]
The last equality is a cyclic reindexing.

Next choose a positive term from the first factor and a negative term from
the second.  Detection means
\[
 2^i-2^{j+r}\equiv-2^k\pmod N,
\]
or equivalently
\[
 2^i+2^k\equiv2^{j+r}\pmod N.
\]
By the uniqueness of a sum of two binary powers, necessarily
\[
 k\equiv i,
 \qquad
 j+r\equiv i+1
 \pmod n.
\]
This contributes
\[
 \sum_i c_id_{i+1-r}.
\]
For the opposite mixed choice, the same argument gives
\[
 \sum_i c_id_{i+r+1}.
\]
Two positive choices do not contribute.  Consequently
\begin{equation}\label{eq:relative-autocorrelation-raw}
 \sum_i d_id_{i+r}
 +\sum_i c_id_{i+1-r}
 +\sum_i c_id_{i+r+1}=0.
\end{equation}

Recall that $b_i=d_i+c_{i-1}$.  Expanding the cyclic autocorrelation of
$b$ gives
\begin{align*}
 \sum_i b_ib_{i+r}
 ={}&\sum_i d_id_{i+r}
 +\sum_i d_ic_{i+r-1}
 +\sum_i c_{i-1}d_{i+r}
 +\sum_i c_{i-1}c_{i+r-1}.
\end{align*}
After reindexing,
\[
 \sum_i d_ic_{i+r-1}=\sum_i c_id_{i+1-r},
 \qquad
 \sum_i c_{i-1}d_{i+r}=\sum_i c_id_{i+r+1},
\]
and
\[
 \sum_i c_{i-1}c_{i+r-1}=\sum_i c_ic_{i+r}.
\]
Therefore \eqref{eq:relative-autocorrelation-raw} is equivalent to
\begin{equation}\label{eq:relative-autocorrelation}
 \sum_i b_ib_{i+r}=\sum_i c_ic_{i+r}
 \qquad(r\in\mathbb Z/n\mathbb Z).
\end{equation}

It remains to encode these identities in the cyclic group algebra.  We have
\begin{align*}
 B(Z)B(Z^{-1})
 &=\left(\sum_i b_iZ^i\right)\left(\sum_jb_jZ^{-j}\right)\\
 &=\sum_{i,j}b_ib_jZ^{i-j}\\
 &=\sum_{r\in\mathbb Z/n\mathbb Z}
   \left(\sum_i b_ib_{i+r}\right)Z^{-r}.
\end{align*}
Indeed, the coefficient of $Z^{-r}$ is obtained from $j=i+r$.  Similarly,
\[
 C(Z)C(Z^{-1})
 =\sum_{r\in\mathbb Z/n\mathbb Z}
   \left(\sum_i c_ic_{i+r}\right)Z^{-r}.
\]
Comparing coefficients and using \eqref{eq:relative-autocorrelation} proves
\eqref{eq:relative-quadratic}.  Conversely, coefficient comparison recovers
all the identities in \eqref{eq:relative-autocorrelation}, so the two
formulations are equivalent.  Finally, evaluation at an $n$th root of unity
gives \eqref{eq:relative-quadratic-eval}.
\end{proof}

\subsection{The relative cubic identity}

The cubic moment supplies the second input.  For an $n$th root of unity $z$, abbreviate
\[
 C_z=C(z),\qquad D_z=D(z),\qquad B_z=B(z).
\]

\begin{proposition}[Relative cubic identity]\label{prop:relative-cubic}
Let $p,q,r\in\overline{\F}_2^*$ satisfy
\[
 p^n=q^n=r^n=1,
 \qquad
 pqr=1.
\]
Then
\begin{align}
0={}&(p^{-1}+q^{-1}+r^{-1})D_pD_qD_r
 +C_pD_qD_r+D_pC_qD_r+D_pD_qC_r\notag\\
&+p^{-1}(p^{-1}+q+r)D_pC_qC_r\notag\\
&+q^{-1}(q^{-1}+r+p)C_pD_qC_r\notag\\
&+r^{-1}(r^{-1}+p+q)C_pC_qD_r.
\label{eq:relative-cubic-D}
\end{align}
Equivalently, after substituting $D_z=B_z+zC_z$,
\begin{align}
0=\Phi_{B,C}(p,q,r):={}&
\alpha\bigl(B_pB_qB_r+C_pC_qC_r\bigr)\notag\\
&+r(p^{-1}+q^{-1})B_pB_qC_r\notag\\
&+q(p^{-1}+r^{-1})B_pC_qB_r\notag\\
&+p(q^{-1}+r^{-1})C_pB_qB_r\notag\\
&+p^{-1}(q+r)B_pC_qC_r\notag\\
&+q^{-1}(p+r)C_pB_qC_r\notag\\
&+r^{-1}(p+q)C_pC_qB_r,
\label{eq:relative-cubic-BC}
\end{align}
where
\[
 \alpha=p^{-1}+q^{-1}+r^{-1}.
\]
\end{proposition}

\begin{proof}
For $u,v\in\mathbb Z/n\mathbb Z$, bijectivity of $f$ gives
\begin{equation}\label{eq:relative-third-moment}
 \sum_{x\in\mathcal T_1}f(x)^{1+2^u+2^v}
 =\sum_{y\in\mathcal T_\varepsilon}y^{1+2^u+2^v}=0.
\end{equation}
The final equality follows from the same trace-fiber power-sum argument as
in the quadratic case.  Since $n\ge5$, the relevant binary sums cannot be
congruent to a residue detected by either trace fiber.  In particular, the
all-positive sign pattern in the expansion below contributes zero.

Write the remaining part of the left-hand side of
\eqref{eq:relative-third-moment} as $E(u,v)$.  Thus
\[
 E(u,v)=0
 \qquad(u,v\in\mathbb Z/n\mathbb Z).
\]
We take its two-dimensional cyclic Fourier transform
\begin{equation}\label{eq:relative-Fourier-transform}
 \widehat E(X,Y)=\sum_{u,v\in\mathbb Z/n\mathbb Z}E(u,v)X^uY^v
\end{equation}
in
\[
 \F_2[X^{\pm1},Y^{\pm1}]/(X^n-1,Y^n-1).
\]
Introduce
\[
 p=XY,
 \qquad
 q=X^{-1},
 \qquad
 r=Y^{-1}.
\]
Then $pqr=1$, and the key separation identity is
\begin{equation}\label{eq:relative-variable-separation}
 p^iq^jr^k=X^{i-j}Y^{i-k}.
\end{equation}
More generally, suppose a carry condition gives
\[
 u=i-j+a,
 \qquad
 v=i-k+b,
\]
where $a,b$ are fixed.  For coefficient sequences
$(\xi_i),(\eta_j),(\zeta_k)$, the corresponding contribution to
\eqref{eq:relative-Fourier-transform} is
\begin{align}
&\sum_{i,j,k}\xi_i\eta_j\zeta_k
 X^{i-j+a}Y^{i-k+b}\notag\\
&\qquad =X^aY^b
 \left(\sum_i\xi_ip^i\right)
 \left(\sum_j\eta_jq^j\right)
 \left(\sum_k\zeta_kr^k\right).
\label{eq:relative-separation-general}
\end{align}
Thus a positive choice at a position produces the corresponding $C$-factor,
whereas a negative choice produces the corresponding $D$-factor.

We now compute the seven non-all-positive sign patterns.

\smallskip
\noindent\emph{Three negative choices.}
The exponent is
\[
 -2^i-2^{j+u}-2^{k+v}.
\]
The three-power carry rule
\[
 2^t+2^t+2^{t+1}\equiv2^{t+2}\pmod N
\]
gives precisely three cyclic possibilities:
\[
\begin{array}{c|c|c|c}
\text{carry condition}&u&v&\text{Fourier factor}\\ \hline
 i=j+u,\ k+v=i+1&i-j&i-k+1&Y=r^{-1}\\
 i=k+v,\ j+u=i+1&i-j+1&i-k&X=q^{-1}\\
 j+u=k+v,\ i=j+u+1&i-j-1&i-k-1&X^{-1}Y^{-1}=p^{-1}.
\end{array}
\]
Using \eqref{eq:relative-separation-general}, the total all-negative
contribution is therefore
\begin{equation}\label{eq:relative-three-negative}
 (p^{-1}+q^{-1}+r^{-1})D(p)D(q)D(r).
\end{equation}

\smallskip
\noindent\emph{One positive and two negative choices.}
We first choose the positive term from the first factor.  The exponent is
\[
 2^i-2^{j+u}-2^{k+v}.
\]
It is detected if and only if, for some $\ell$,
\[
 2^i+2^\ell=2^{j+u}+2^{k+v}.
\]
By uniqueness of a sum of two binary powers, this is equivalent to
\[
 j+u=i
 \quad\text{or}\quad
 k+v=i.
\]
Over $\F_2$, the indicator of this union is
\[
 \delta_{j+u,i}+\delta_{k+v,i}
 +\delta_{j+u,i}\delta_{k+v,i}.
\]
For the first single condition, the index $k$ is free, so its total contains
$\sum_kd_k=D(1)=0$.  The second single condition vanishes in the same way.
Only their intersection remains; it gives
\[
 u=i-j,
 \qquad
 v=i-k,
\]
and hence contributes
\[
 C(p)D(q)D(r).
\]

For completeness, consider the middle sign pattern $(-+-)$ separately.
Its exponent is
\[
 -2^i+2^{j+u}-2^{k+v},
\]
and detection is equivalent to
\[
 2^{j+u}+2^\ell=2^i+2^{k+v}.
\]
Thus
\[
 j+u=i
 \quad\text{or}\quad
 j+u=k+v.
\]
The contribution of the first condition vanishes after summing over the
free index $k$, and that of the second vanishes after summing over the free
index $i$; both contain $D(1)$.  Their intersection gives
$u=i-j$ and $v=i-k$, and therefore contributes
\[
 D(p)C(q)D(r).
\]
The third placement is analogous.  Hence the total contribution of the
three patterns with one positive sign is
\begin{equation}\label{eq:relative-one-positive}
 C(p)D(q)D(r)
 +D(p)C(q)D(r)
 +D(p)D(q)C(r).
\end{equation}

\smallskip
\noindent\emph{One negative and two positive choices.}
First place the negative term in the first factor.  The exponent is
\[
 -2^i+2^{j+u}+2^{k+v}.
\]
Detection means that, for some $\ell$,
\[
 2^{j+u}+2^{k+v}+2^\ell=2^i.
\]
The three possible carry patterns are
\[
\begin{array}{c|c|c|c}
\text{carry condition}&u&v&\text{Fourier factor}\\ \hline
 j+u=k+v=i-2&i-j-2&i-k-2&X^{-2}Y^{-2}=p^{-2}\\
 j+u=i-2,\ k+v=i-1&i-j-2&i-k-1&X^{-2}Y^{-1}=p^{-1}q\\
 j+u=i-1,\ k+v=i-2&i-j-1&i-k-2&X^{-1}Y^{-2}=p^{-1}r.
\end{array}
\]
Consequently this sign pattern contributes
\begin{align*}
 &(p^{-2}+p^{-1}q+p^{-1}r)
 D(p)C(q)C(r)\\
 &\qquad =p^{-1}(p^{-1}+q+r)D(p)C(q)C(r).
\end{align*}
The other two placements are obtained cyclically.  Adding
\eqref{eq:relative-three-negative},
\eqref{eq:relative-one-positive}, and these three contributions proves
\eqref{eq:relative-cubic-D} for triples of the form
$(p,q,r)=(XY,X^{-1},Y^{-1})$.

It remains to verify that this parametrization covers every triple
appearing in the statement.  Let $p,q,r$ be arbitrary $n$th roots of unity
satisfying $pqr=1$.  Setting $X=q^{-1}$ and $Y=r^{-1}$, we have
$X^n=Y^n=1$ and
\[
 XY=(qr)^{-1}=p,\qquad X^{-1}=q,\qquad Y^{-1}=r.
\]
Thus every admissible triple is of the form $(XY,X^{-1},Y^{-1})$.
Consequently, the identity derived above holds for all $p,q,r$ satisfying
the hypotheses of the proposition.

It remains to derive the $B,C$ form.  Set
\[
 \alpha=p^{-1}+q^{-1}+r^{-1}
\]
and
\[
 \lambda_p=p^{-1}(p^{-1}+q+r),
 \qquad
 \lambda_q=q^{-1}(q^{-1}+r+p),
 \qquad
 \lambda_r=r^{-1}(r^{-1}+p+q).
\]
Since $B(Z)=D(Z)+ZC(Z)$ and the characteristic is two,
\[
 D_p=B_p+pC_p,
 \qquad
 D_q=B_q+qC_q,
 \qquad
 D_r=B_r+rC_r.
\]
Substitute these expressions into \eqref{eq:relative-cubic-D}.  We collect
the coefficients of the eight possible monomial types.

The coefficient of $B_pB_qB_r$ is $\alpha$.  The coefficient of
$B_pB_qC_r$ is
\[
 \alpha r+1=r(p^{-1}+q^{-1}),
\]
and cyclic permutation gives the other two terms with two $B$-factors.

For $B_pC_qC_r$, the four sources are
\[
 \alpha qr,
 \qquad r,
 \qquad q,
 \qquad \lambda_p.
\]
Since $pqr=1$,
\[
 \alpha qr=p^{-2}+r+q,
\]
and therefore
\begin{align*}
 \alpha qr+q+r+\lambda_p
 &=p^{-2}+p^{-1}(p^{-1}+q+r)\\
 &=p^{-1}(q+r).
\end{align*}
Cyclic permutation gives the other two terms with one $B$-factor.

Finally, consider $C_pC_qC_r$.  The first four terms of
\eqref{eq:relative-cubic-D} contribute
\[
 \alpha pqr+qr+pr+pq
 =\alpha+(p^{-1}+q^{-1}+r^{-1})=0.
\]
The last three terms contribute
\[
 p\lambda_p+q\lambda_q+r\lambda_r.
\]
Here the correct identities are
\[
 p\lambda_p=p^{-1}+q+r,
 \qquad
 q\lambda_q=q^{-1}+r+p,
 \qquad
 r\lambda_r=r^{-1}+p+q.
\]
The terms $p,q,r$ each occur twice and cancel, leaving
\[
 p\lambda_p+q\lambda_q+r\lambda_r
 =p^{-1}+q^{-1}+r^{-1}=\alpha.
\]
Thus the coefficient of $C_pC_qC_r$ is also $\alpha$.  Collecting all eight
monomial types gives exactly \eqref{eq:relative-cubic-BC}.
\end{proof}

\subsection{Spectral-zero rigidity}

We first exclude a common nontrivial zero of $B$ and $C$.

\begin{lemma}[No common spectral zero]\label{lem:no-common-zero}
If $z^n=1$ and $z\ne1$, then
\[
 (B(z),C(z))\ne(0,0).
\]
\end{lemma}

\begin{proof}
Suppose that $B(z)=C(z)=0$, and let $h(Z)$ be the minimal polynomial of $z$ over $\F_2$.
Then
\[
 h\mid Z^n-1,
 \qquad
 h\ne Z+1.
\]
Since $B(z)=D(z)+zC(z)$, we also have $D(z)=0$.  Thus
\[
 C=hC_1,
 \qquad
 D=hD_1
\]
for some $C_1,D_1\in\F_2[Z]$.  Because $h(1)=1$ and $D(1)=0$, one has $D_1(1)=0$.
Define
\[
 g(x)=L_{C_1}(x)+L_{D_1}(x^{-1}).
\]
Composition of binary linearized polynomials corresponds to multiplication of their coefficient
polynomials modulo $Z^n-1$, and here it gives
\begin{equation}\label{eq:factor-through-h}
 f=L_h\circ g.
\end{equation}

The injectivity of $f$ on $\mathcal T_1$ implies the injectivity of $g$ there.  Since
$D_1(1)=0$, the image $g(\mathcal T_1)$ is one trace fiber, say $\mathcal T_\eta$; equality
of cardinalities shows that
\[
 g(\mathcal T_1)=\mathcal T_\eta.
\]
Equation \eqref{eq:factor-through-h} then forces $L_h$ to be injective on
$\mathcal T_\eta$, and hence on both trace fibers by
\Cref{lem:trace-fiber-kernel}.  But $h\mid Z^n-1$ implies that $L_h$ has nonzero kernel.
Moreover, $h(1)=1$, so \Cref{lem:trace-fiber-kernel} says that trace-fiber injectivity would
force $L_h$ to be a permutation of $F$.  This contradiction proves the lemma.
\end{proof}

The next lemma is the central local consequence of the two relative identities.

\begin{lemma}[Two-zero rigidity]\label{lem:two-zero-rigidity}
Let $p^n=q^n=1$ and suppose
\[
 p\ne1,
 \qquad
 q\ne1,
 \qquad
 pq\ne1.
\]
If $C(p)=C(q)=0$, then $p=q$.
\end{lemma}

\begin{proof}
Put
\[
 r=(pq)^{-1}.
\]
By \Cref{lem:no-common-zero},
\[
 B(p)B(q)\ne0.
\]
The quadratic identity \eqref{eq:relative-quadratic-eval} gives
\[
 B(p^{-1})=B(q^{-1})=0,
\]
and \Cref{lem:no-common-zero} then gives
\begin{equation}\label{eq:inverse-C-nonzero}
 C(p^{-1})C(q^{-1})\ne0.
\end{equation}
Set
\[
 \alpha=p^{-1}+q^{-1}+r^{-1},
 \qquad
 \alpha'=p+q+r,
 \qquad
 \beta=r(p^{-1}+q^{-1}).
\]
Substituting $C(p)=C(q)=0$ into \eqref{eq:relative-cubic-BC} leaves
\[
 B(p)B(q)\bigl(\alpha B(r)+\beta C(r)\bigr)=0.
\]
Hence
\begin{equation}\label{eq:alpha-beta-relation}
 \alpha B(r)+\beta C(r)=0.
\end{equation}

If $B(r)=0$, then \Cref{lem:no-common-zero} gives $C(r)\ne0$, so
$\beta=0$.  Thus $p^{-1}+q^{-1}=0$, and hence $p=q$.

Suppose next that $C(r)=0$.  Then $B(r)\ne0$ and
\eqref{eq:alpha-beta-relation} gives $\alpha=0$.  The quadratic identity gives
$B(r^{-1})=0$, and the no-common-zero lemma gives $C(r^{-1})\ne0$.  Apply
\eqref{eq:relative-cubic-BC} to the triple
\[
 (p^{-1},r^{-1},q^{-1}).
\]
At all three positions $B$ vanishes and $C$ does not, so the identity reduces to
$\alpha'=0$.  However,
\begin{equation}\label{eq:alpha-factorization}
 pq(\alpha+\alpha')
 =(p+1)(q+1)(pq+1),
\end{equation}
which is nonzero under the hypotheses.  This is a contradiction.

It remains to consider $B(r)C(r)\ne0$.  Put
\[
 u=\frac{B(r)}{C(r)}.
\]
Then \eqref{eq:alpha-beta-relation} becomes
\begin{equation}\label{eq:u-first}
 \alpha u+\beta=0.
\end{equation}
By the quadratic identity,
\[
 \frac{B(r^{-1})}{C(r^{-1})}=u^{-1}.
\]
Applying the cubic identity to $(p^{-1},r^{-1},q^{-1})$ and using
\eqref{eq:inverse-C-nonzero} gives
\[
 \alpha'+\beta u^{-1}=0,
\]
or
\begin{equation}\label{eq:u-second}
 \alpha'u+\beta=0.
\end{equation}
Adding \eqref{eq:u-first} and \eqref{eq:u-second} yields
$(\alpha+\alpha')u=0$.  Since $u\ne0$, this contradicts
\eqref{eq:alpha-factorization}.  The only possible case is therefore the first one, and it
gives $p=q$.
\end{proof}

\begin{proposition}\label{prop:no-nontrivial-C-zero}
If $z^n=1$ and $z\ne1$, then $C(z)\ne0$.
\end{proposition}

\begin{proof}
Suppose that $C(z)=0$.  Since $C$ has coefficients in $\F_2$,
\[
 C(z^2)=C(z)^2=0.
\]
If $z^3\ne1$, apply \Cref{lem:two-zero-rigidity} with $p=z$ and $q=z^2$.  Its
hypotheses hold, and it gives $z=z^2$, hence $z=1$, a contradiction.

If $z^3=1$, then $z\ne1$ implies $z^{-1}=z^2$.  Thus
\[
 C(z)=C(z^{-1})=0.
\]
By \Cref{lem:no-common-zero}, both $B(z)$ and $B(z^{-1})$ are nonzero, whereas the
quadratic identity gives
\[
 B(z)B(z^{-1})=C(z)C(z^{-1})=0.
\]
This is again a contradiction.
\end{proof}

\subsection{Proof of Conjecture 4.5 and the structural classification}

\begin{proof}[Proof of \Cref{thm:general-reduction}]
By \Cref{prop:no-nontrivial-C-zero}, the polynomial $C$ has no zero among the nontrivial
$n$th roots of unity.

If $C(1)=1$, then $1$ is not a zero either, and therefore
\[
 \gcd(C(Z),Z^n-1)=1.
\]
By \Cref{lem:trace-fiber-kernel}, $L_C\in\mathscr L_n$.

Suppose that $C(1)=0$.  By \Cref{lem:even-degree-parity}, $n$ is odd.  Hence
$Z^n-1$ is squarefree in characteristic two.  The only common root of $C(Z)$ and
$Z^n-1$ is then $1$, and it occurs with multiplicity one in $Z^n-1$.  Consequently
\[
 \gcd(C(Z),Z^n-1)=Z+1.
\]
The second case of \Cref{lem:trace-fiber-kernel} again gives
$L_C\in\mathscr L_n$.  This proves Conjecture~4.5.
\end{proof}

Let $*$ denote cyclic binary convolution of coefficient words, so that
\[
 L_a\circ L_b=L_{a*b}.
\]
The following consequence combines the two Hollmann--Xiang conjectures.

\begin{corollary}[All Kloosterman polynomials]\label{cor:all-kloosterman}
Let $n>3$.  A function $L_{c,d}$ is a Kloosterman polynomial on $\F_{2^n}$ if and only if
there exist $L_a\in\mathscr L_n$ and
\[
 e\in\{0,3,6,10\}
\]
such that
\[
 L_{c,d}=L_a\circ L_{1,e}.
\]
\end{corollary}

\begin{proof}
Suppose first that $L_{c,d}$ is Kloosterman.  By
\Cref{thm:general-reduction}, $L_c\in\mathscr L_n$.  The set $\mathscr L_n$ is a group
under composition by \cite[Lemma~4.3]{HollmannXiang2004}; choose $L_b\in\mathscr L_n$
with
\[
 L_b\circ L_c=L_1.
\]
Then
\[
 L_b\circ L_{c,d}=L_{1,b*d}
\]
is again Kloosterman.  By \Cref{thm:main},
\[
 b*d\in\{0,3,6,10\}.
\]
Writing $e=b*d$ and composing with $L_c$ gives
\[
 L_{c,d}=L_c\circ L_{1,e}.
\]

Conversely, $L_{1,e}$ is Kloosterman for each listed $e$ by
\cite[Theorem~4.1]{HollmannXiang2004} and the trivial case $e=0$.  Composing with any
$L_a\in\mathscr L_n$ preserves the Kloosterman property by
\cite[Lemma~4.3(ii)]{HollmannXiang2004}.
\end{proof}

\section{Concluding remarks}

We have proved both conjectures proposed by Hollmann and Xiang for binary Kloosterman
polynomials.  For the normalized family $L_{1,u}$, their sufficiency theorem and the necessity
proved here give the exact list
\[
 u\in\{0,3,6,10\}.
\]
The necessity is governed by a rigid sequence of moment constraints: the quadratic moment
gives a cyclic inverse identity, the cubic moment reduces the coefficient word to a monomial
or its complement, and a small number of explicit cubic and quartic moments completes the
elimination.

For a general Kloosterman polynomial $L_{c,d}$, the relative quadratic identity
\[
 B(Z)B(Z^{-1})=C(Z)C(Z^{-1}),
 \qquad B(Z)=D(Z)+ZC(Z),
\]
and the relative cubic spectral identity force $C$ to have no nontrivial $n$th root of unity as
a zero.  This proves that $L_c$ is injective on both trace fibers and establishes Conjecture~4.5.
Together with the normalized classification, it follows that every Kloosterman polynomial is
a composition
\[
 L_a\circ L_{1,e},
 \qquad
 L_a\in\mathscr L_n,
 \qquad
 e\in\{0,3,6,10\}.
\]
Thus the two conjectures provide not only a list of the basic examples but also a structural
classification of all binary Kloosterman polynomials.

The relative-moment method suggests a broader problem: for rational functions built from
several linearized pieces, one may seek quadratic identities comparing their spectral norms and
higher-moment identities that rigidify the possible zero sets.  Similar questions on other trace
fibers or in odd characteristic would require different power-sum supports and carry rules, but
the same combination of cyclic group algebra and spectral-zero analysis may remain useful.

\end{document}